\documentclass{amsart}

\usepackage{amssymb}
\usepackage{mathrsfs, textcomp}
\usepackage{graphicx}
\usepackage{bm}
\usepackage{enumitem}

\usepackage{amsmath,amsfonts,amsthm,latexsym}
\usepackage[all]{xy}

\newtheorem{theorem}{Theorem}[section]
\newtheorem{lemma}[theorem]{Lemma}
\newtheorem{proposition}[theorem]{Proposition}

\theoremstyle{definition}
\newtheorem{definition}[theorem]{Definition}
\newtheorem{example}[theorem]{Example}

\numberwithin{equation}{section}

\newcommand{\R}{\mathbb{R}}

\newcommand{\cC}{\mathcal{C}}

\newcommand{\cF}{\mathcal{F}}

\newcommand{\cP}{\mathcal{P}}
\newcommand{\cR}{\mathcal{R}}

\newcommand{\cV}{\mathcal{V}}

\newcommand{\fB}{{\mathfrak B}}

\newcommand{\sG}{\mathscr{G}}

\newcommand{\ssm}{\smallsetminus}
\newcommand{\pperp}{\perp\perp}
\newcommand{\Int}{{\rm int \ }}

\DeclareMathOperator{\Min}{Min}

\DeclareMathOperator{\Max}{Max}

\DeclareMathOperator{\cl}{cl}
\DeclareMathOperator{\coz}{coz}

\DeclareMathOperator{\Spec}{Spec}

\DeclareMathOperator{\Ult}{{\bf Ult}}
\DeclareMathOperator{\Coz}{Coz}

\newcommand{\ra}{\rightarrow}

\theoremstyle{definition}
\theoremstyle{definition}
\theoremstyle{definition}
\theoremstyle{definition}
\theoremstyle{definition}
\theoremstyle{definition}
\theoremstyle{remark}
\theoremstyle{definition}
\theoremstyle{definition}

\begin{document}

\title{When $\Max_d(G)$ is zero-dimensional }
\author[R. Carrera]{Ricardo Carrera${}^1$}
\address{${}^1$Department of Mathematics, Halmos College of Arts and Sciences, Nova Southeastern University, Fort Lauderdale, FL 33328}
\email{ricardo@nova.edu}
\author[R.~Lafuente-Rodriguez]{Ramiro Lafuente-Rodriguez${}^2$}
\address{${}^2$University of South Dakota, Department of Mathematical Sciences, Vermillion, SD 57069}
\email{Ramiro.LafuenteRodri@usd.edu}
\author[W. Wm. McGovern]{Warren Wm.~McGovern${}^3$}
\address{${}^3$H. L. Wilkes Honors College, Florida Atlantic University, Jupiter, FL 33458}
\email{warren.mcgovern@fau.edu}
\keywords{archimedean lattice-ordered group, $d$-subgroups; zero-dimensional compact Hausdorff space}

\begin{abstract}
This article is a continuation of \cite{lm} where a classification of when the space of minimal prime subgroups of a given lattice-ordered group equipped with the inverse topology has a clopen $\pi$-base. For nice $\ell$-groups, (e.g. {\bf W}-objects) this occurs precisely when the space of maximal $d$-subgroups ({\it qua} the hull kernel topology) has a clopen $\pi$-base. It occurred to us that presently there is no classification of when the space of maximal $d$-subgroups of a {\bf W}-object is zero-dimensional, except for the case of the $C(X)$, the real-valued continuous functions on a topological space $X$, considered in \cite{hvw2}. We generalize the situation without the aid of the cloz-cover of a compact Hausdorff space.

\end{abstract}

\maketitle

\thispagestyle{empty}

\section{Introduction and Preliminaries}

Recall that {\bf W} is the category whose objects are pairs $(G,u)$ such that $G$ is an archimedean $\ell$-group and $u\in G$ is a positive weak order unit. We assume the reader is familiar with the theory of lattice-ordered groups, and in particular {\bf W}-objects. We use $\cC(G)$ to denote the collection of convex $\ell$-subgroups of $G$. This is an algebraic frame with FIP. In $G$, there are convex $\ell$-subgroups that are maximal with respect to not containing $u$: Zorn's Lemma. Such subgroups are called {\bf values of $u$}, and the set of these is called the {\bf Yosida space of $(G,u)$}; it is denoted by $YG$. When equipped with the hull-kernel topology,  $YG$ is a compact Hausdorff space. The collection of sets of the form
$$\coz(g)=\{V\in YG: g\notin V\}$$
for some $g\in G$, forms a base for the open subsets; any set of the form $\coz(g)$ is called a $G$-cozeroset of $YG$. The complement of $\coz(g)$ is denoted by $Z(g)$ and is called the {\bf zero-set of $g$}. The collections of $G$-cozerosets of $YG$ and $G$-zero-sets of $G$ will be denoted by $\Coz(G)$ and $Z(G)$, respectively. The collection of $G$-cozerosets is closed under finite union and finite intersection. Namely, for $g,h\in G^+$, $\coz(g)\cap \coz(h)=\coz(g\wedge h)$ and $\coz(g)\cup\coz(h)=\coz(g+h)=\coz(g\vee h)$.

The main tool in analysing {\bf W}-objects is the Yosida representation. Let $\overline{R}=\R\cup\{\pm\infty\}$, the two-point compactification of the reals. For a compact Haudorff space $X$, a continuous function $f:X\ra \overline{R}$ is said to be {\bf almost real-valued} if $f^{-1}(\R)$ is a dense subset of $X$. The collection of all almost real-valued continuous functions on $X$ is denoted by $D(X)$. In general, $D(X)$ need not be a group as addition might not be well-defined. However, we do speak of $\ell$-subgroups of $D(X)$. For example, the collection of continuous real-valued functions defined on $X$ is an $\ell$-subgroup of $D(X)$.

\begin{theorem}[The Yosida Representation]
Let $(G,u)\in {\bf W}$. There is an $\ell$-group $\hat{G}$ of continuous almost real-valued function on $YG$ and an $\ell$-isomorphism of $G$ onto $\hat{G}$ such that $u\mapsto {\bf 1}$. Furthermore, for each closed subset $V$ of $YG$ and $p\notin V$ there is some $g\in G$ such that $\hat{g}(p)=1$ and $\hat{g}(V)=0$. $YG$ is the unique (up to homeomorphism) compact Hausdorff space with this property.
\end{theorem}

We identify any $(G,u)$ with $(\hat{G},{\bf 1})$. The property involving points not in closed sets will be referred to as the separation property.
\vspace{.2in}

We view $C(X)$ as a {\bf W}-object by taking the constant function {\bf 1} as its distinguished weak order unit: $(C(X),{\bf 1})$. Recall that if $X$ is a compact Hausdorff space then the Yosida space of $X$ is in fact homeomorphic to $X$ and the values of $(C(X),{\bf 1})$ are precisely sets of the form
$M_p=\{f\in C(X):f(p)=0\}$
for some $p\in X$. In this case the zero-sets of $(C(X),{\bf 1})$ are precisely the usual zero-sets of $X$. Because of this we typically view the values of $u$ as $V_p$ for some $p\in YG$; namely those elements $g\in G$ for which $p$ is in the $G$-zeroset of $g$.

Values are examples of {\bf prime} subgroups: a convex $\ell$-subgroup $P$ with the property that if $a,b\in G^+$ and $a\wedge b=0$, then either $a\in P$ or $b\in P$. Primes can be characterized as those $P$ whose factor set of cosets are totally ordered. The collection of all prime subgroups is known as the prime spectrum of $G$ and is denoted by $\Spec(G)$. Lattice-ordered groups have the feature that the collection of primes containing a given prime form a chain, i.e. $\Spec(G)$ is a root system.

The collection of prime subgroups is also equipped with the hull-kernel topology: the collection of sets of the form $U(g)=\{P\in\Spec(G): g\notin P\}$ forms a base for the open sets. The collection of minimal prime subgroups of $G$ is denoted by $\Min(G)$. Investigating topological properties of $\Min(G)$ has been an important part of the theory of ${\bf W}$-objects. For example, $\Min(G)$ can also be equipped with the inverse topology which is the topology generated by the collection of sets of the form $V(g)=\Min(G)\ssm U(g)$; we use $\Min(G)^{-1}$ to denote the space of minimal prime subgroups equipped with the inverse topology. In some sense, this article is about classifying spaces of prime subgroups.

Next, given $S\subseteq G$, the polar of $S$ is denoted by $S^\perp$ and defined as
$$S^\perp = \{h\in G: |h|\wedge |g|=0 \hbox{ for all } h\in S\}.$$
A polar subgroup  is a subgroup of the form $S^\perp$ for some subset $S\subseteq G$. When $S=\{g\}$ we instead write $g^\perp$. A {\bf principal polar} is a set of the form $g^{\pperp}$ for some $g\in G$. If $g^\perp=0$, then $g$ is called a {\bf weak order unit}. The set of polars of $G$ is denoted by $\cP(G)$, and when ordered by inclusion it is a complete boolean algebra; more on this later.

A convex $\ell$-subgroup $H$ of $G$ is called a {\bf $d$-subgroup} if whenever $h\in H$, then $h^{\pperp}\subseteq H$. No proper $d$-subgroup can contain a weak order unit, thus maximal $d$-subgroups exist; denote the set of all maximal $d$-subgroups by $\Max_d(G)$. Proposition 4.3 of \cite{bm4} classifies maximal $d$-subgroups as being those convex $\ell$-subgroups which are maximal with respect to not containing any weak order units.

As we do for the Yosida space, we equip $\Max_d(G)$ with the hull-kernel topology. For $g\in G$, let $U_d(g)=\{M\in\Max_d(G): g\notin M\}$. Clearly, $U_d(g)=U_d(|g|)$. The set-theoretic complement of $U_d(g)$ in $\Max_d(G)$ shall be denoted by $V_d(g)$.
\bigskip

\begin{lemma}\cite[Lemma 4.4]{bm4}
Let $(G,u)\in {\bf W}$.  The following hold for all $g,h\in G^+$.
\begin{enumerate}[label={\rm (\alph*)}, nolistsep]
\item $U_d(g)=\Max_d(G)$ if and only if $g$ is a weak order unit.
\item $U_d(g)=\emptyset$ if and only if $g=0$.
\item $U_d(g)\cup U_d(h)=U_d(g\vee h)$.
\item $U_d(g)\cap U_d(h)=U_d(g\wedge h)$.
\item The subset $T\subseteq \Max_d(G)$ is open in the hull kernel topology if and only if there is some $d$-subgroup $H$ for which $T=U_d(H)$.
\end{enumerate}
\end{lemma}

In \cite{lm}, the authors showed that if $G$ is a Lamron $\ell$-group, then $\Min(G)^{-1}$ is zero-dimensional if and only if $\Max_d(G)$ is zero-dimensional. We were able to show that in the case of ${\bf W}$-objects, $\Min(G)^{-1}$ has a clopen $\pi$-base if and only if $\Max_d(G)$ has a clopen $\pi$-base. In what follows, our work is aimed at answering when $\Max_d(G)$ is zero-dimensional.
\vspace{.2in}

In any lattice-ordered group, $\Spec(G)$ is a root system. In the case of a {\bf W}-object, say $(G,u)$, for each $M\in \Max_d(G)$ we know that $u\notin M$ and so $M$ can be extended to a unique value of $u$. This defines a map
$$\lambda_u:\Max_d(G)\ra YG$$
which is continuous.

In what follows we shall need to use properties of the operator $U_d(\cdot)$. Let $X$ denote a compact Hausdorff space and let $\cR(X)$ denote the collections of regular closed subsets of $X$. (Recall that $V\subseteq X$ is called {\it regular closed} if $V=\cl_X \Int _X V$.) It is well-known that $\cR(X)$ is a (complete) boolean algebra when partially ordered by inclusion. The meet, join, and complement are given as follows: for $V_1,V_2\in \cR(X)$
\medskip

\begin{enumerate}[label={\rm (\roman*)}, nolistsep]
\item $V_1\cup' V_2 = V_1\cup V_2$;
\item $V_1\cap' V_2 = \cl_X \Int _X (V_1\cap V_2)$;
\item $V_1'= \cl_X (X\ssm V_1)$.
\end{enumerate}

Certain bounded sublattices of $\cR(X)$ are of importance. We recall the definition of a Wallman lattice (see \cite[Definition 2.1]{bm4}). The condition (vi) was not included in the article, but we need it here. Wallman described this property as the {\it disjunction} property in \cite[Lemma 3.]{wallman}.

\begin{definition}
Let $(L,\vee,\wedge,0,1)$ be a bounded distributive lattice, and let $\Ult(L)$ denote the collection of $L$-ultrafilters. For $a\in L$, denote the set of ultrafilters containing $a$ by  $\cV(a)$. The operator $\cV(\cdot)$ has the following properties: \cite{banaschewski}.

\begin{enumerate}[label={\rm (\roman*)}, nolistsep]\label{wallman}
\item For each $a,b\in L$, $\cV(a)\cup \cV(b)=\cV(a\vee b)$ and $\cV(a)\cap \cV(b)=\cV(a\wedge b)$.
\item The collection $\{\cV(a):a\in L\}$ forms a base for a topology of closed sets on $\Ult(L)$. This is called the {\bf Wallman topology} on $\Ult(L)$.
\item For each $a<1$, there is a $0<c\in L$ such that $a\wedge c=0$ if and only if the map $a\ra \cV(a)$ is injective. (A lattice satisfying either of these equivalent conditions is called a {\bf Wallman lattice}.)
\item If $L$ is a Wallman lattice, then $\Ult(L)$ is a compact $T_1$-space.
\item The space $\Ult(L)$ is a Hausdorff space if and only if for any $a,b\in L$ such that $a\wedge b=0$ there exists $x,y\in L$ such that $x\vee y=1$ and $a\wedge y=0=b\wedge x$. (When $\Ult(L)$ is a Hausdorff space, we shall say $L$ is a {\bf normal} lattice.)
\item $L$ is a {\bf disjunctive lattice} if whenever $a\neq b$ there is some $c\in L$ such that either $a\wedge  c=0< b\wedge c$ or $b\wedge c=0<a\wedge c$.
\end{enumerate}
\end{definition}

When we deal with a {\bf W}-object, say $G$, then we are interested in the following sublattices of $\cR(YG)$. (In our view, these are the appropriate generalizations of the objects studied in \cite{hvw2}.)

\label{wl2}
\begin{enumerate}[label={\rm (\alph*)}, nolistsep]
  \item $Z^\sharp(G)=\{\cl \Int Z(f): f\in G\}.$
  \item $\cl\Coz(G) = \{\cl C: C\in \Coz(G)\}.$
  \item $\fB(G)=\{C\subseteq YG : C$ is clopen$\}$.
\end{enumerate}
\vspace{.2in}

The collection $Z^\sharp(G)$ is a normal Wallman sublattice of $\cR(YG)$. In particular, for any pair of $G$-zero-sets $Z_1,Z_2\in Z(G)$.
$$\cl \Int Z_1 \cap' \cl \Int Z_2 = \cl \Int (Z_1 \cap Z_2)$$
and
$$\cl \Int Z_1 \cup' \cl \Int Z_2 = \cl \Int (Z_1 \cup Z_2).$$

Moreover, the collection $Z^\sharp(G)$ is disjunctive.
\begin{lemma}
Let $(G,u)\in {\bf W}$. The lattice $Z^\sharp(G)$ is a disjunctive lattice.
\end{lemma}

\begin{proof}
Let $e,f \in G^+$ such that $\cl\Int Z(e)\neq \cl\Int Z(f)$. Without loss of generality there is some $p\in \cl\Int Z(f) \ssm \cl\Int Z(e)$. Then, by the separation property, there is some $g\in G^+$ such that $\cl\Int Z(e)\subseteq Z(g)$ and $p\notin Z(g)$. Then $\cl\Int (Z(e)\cap Z(g))=\cl\Int Z(e)$ so that by replacing $e$ with $e+g$ we may assume that $p\notin Z(e)$, and then (by using the separation property again) one can find an $h\in G^+$ such that $p\in \Int Z(h)$ and $Z(h)\cap Z(e)=\emptyset$. We leave it to the interested reader to check that $\cl\Int Z(h)$ is disjoint from $\cl\Int Z(e)$ while $\cl\Int Z(h)\cap' \cl\Int Z(f)\neq \emptyset$.
\end{proof}

In \cite[Theorem 4.8]{bm4}, the authors showed that the space of ultrafilters of $Z^\sharp(G)$ is homeomorphic to $\Max_d(G)$. In particular, they showed that for any $d$-subgroup, $H$, the collection $Z^\sharp[H]=\{\cl\Int Z(h): h\in H\}$ is a $Z^\sharp(G)$-filter. Conversely, if $\cF$ is a $Z^\sharp$-filter, then $\overleftarrow{Z}^{\sharp}[\cF]=\{g\in G: \cl\Int Z(g)\in \cF\}$ is a $d$-subgroup of $G$. We can use this for a description of $\lambda_u(M)$.

\begin{lemma}\label{lambda}
Let $M\in \Max_d(G)$. Then $\lambda_u(M)=V_p$ where $\{p\}= \cap \{\cl\Int Z(g): g\in M\}$.
\end{lemma}

\begin{proof}
Let $p\in YG$ such that $\lambda_u(M)=V_p$. Then for each $g\in M$, $g\in V_p$ and $p\in Z(g)$. If for some $m\in M$, $p\notin \cl\Int Z(m)$, then there is some $h\in G^+$ such that $h(p)=1$ and $\cl\Int Z(m) \subseteq \Int Z(h)$. It follows that $\cl\Int Z(m)\subseteq \cl\Int Z(h)$, whence $h\in M$. But by design $h\notin V_p$, the desired contradiction.
\end{proof}

\begin{proposition}
Let $(G,u)\in {\bf W}$ and $g,h\in G^+$. The following are equivalent.
\begin{enumerate}[label={\rm (\arabic*)}, nolistsep]
  \item $U_d(g)=U_d(h)$,
  \item $V_d(g)=V_d(h)$,
  \item $\cl\coz(g)=\cl\coz(h)$,
  \item $\cl\Int Z(g)=\cl\Int Z(h)$,
  \item $g^{\pperp}=h^{\pperp}$.

\end{enumerate}

\end{proposition}

\begin{proof}
That (1) and (2) are equivalent is purely set-theoretic; same for (3) and (4).
That (4) and (5) are equivalent follows from the characterization of $g^{\pperp}$ given in
\cite[Lemma 4.6. (b)]{bm4}, namely
$$g^{\pperp}= \{k\in G: \cl\Int Z(g)\subseteq \cl\Int Z(k)\}.$$
\vspace{.1in}

That (5) implies (1) is obvious using the definition of $d$-subgroup.
\vspace{.1in}

So assume (1). If $\cl\Int Z(g)\neq \cl\Int Z(h)$, then the conjunctive property of $Z^\sharp(G)$ produces an element $f\in G^+$ such that either $\cl \Int Z(g) \cap' \cl\Int Z(f)=\emptyset$ or $\cl \Int Z(h) \cap' \cl\Int Z(f)=\emptyset$. Either way, there must be a maximal $d$-subgroup containing either $g$ and not $h$, or the other way around. This contradicts (1).
\end{proof}

\vspace{.2in}

\begin{definition}
For any $g\in G^+$ , we say that $g$ is {\bf complemented} if there is an $f\in G^+$ such that $g\wedge f=0$ and $g\vee f$ is a weak order unit. We let $c(G)$ denote the set of all complemented element of $G$. The cozeroset of such an element is a {\it complemented cozero} of $YG$, that is, there is another cozeroset, say $D\in \Coz(G)$, for which $C\cap D=\emptyset$ and $C\cup D$ is dense in $X$. The set $c(G)$ is a sublattice of $G^+$. ($G$ is called {\bf complemented} when $G^+=c(G)$. This is a well-known class of $\ell$-groups.)

This leads us to recall the definition of the following sublattice of $Z^\sharp(G)$.
$$\sG(G)=\{\cl \coz(e): e\in c(G)\}.$$

One ought to notice that $\sG(G)=Z^\sharp(G)\cap \cl\coz(G)$ and therfore $\sG(G)$ is a boolean algebra. This lattice will play a role in the next section. We end this section relating the complemented elements of $G$ to the clopen subsets of $\Max_d(G)$.
\end{definition}

\begin{lemma}\cite[Lemma 5.1]{lm}\label{maxclopen}
The set $K\subseteq \Max_d(G)$ is clopen if and only if $K=U_d(e)$ for some $e\in c(G)$.
\end{lemma}
\vspace{.2in}

\begin{proposition}
Let $e\in c(G)$. Then $\lambda_u(V_d(e)) = \cl\Int Z(e)$.
\end{proposition}

\begin{proof}
Let $e\in c(G)$ and choose any $f\in c(G)$ such that $e\wedge f=0$ and $e\vee f$ is a weak order unit.

Let $M\in V_d(e)$ and let $\lambda_u(M)=V_p$. Then by Lemma \ref{lambda}, $p\in \cl\Int Z(e)$, whence $\lambda_u(V_d(e))\subseteq \cl\Int Z(e)$.

Next, we show that $\Int Z(e)\subseteq \lambda_u(V_d(e))$. To that end, let $p\in \Int Z(e)$ and let $M\in \Max_d(G)$ such that $\lambda_u(M)=V_p$. We claim that $e\in M$, i.e. $M\in V_d(e)$. If not, then $f\in M$. Applying Lemma \ref{lambda}, we gather that $p\in \cl\Int Z(f)$, which is a contradiction.

Putting both pieces together produces the following string of containments.
$$\Int Z(e)\subseteq \lambda_u(V_d(e)) \subseteq \cl\Int Z(e).$$
Now, since $\Max_d(G)$ is compact and $V_d(e)$ is clopen it follows that both $V_d(e)$ and $\lambda_u(V_d(e))$ are compact ($\lambda_u$ is continuous). Consequently, $\lambda_u(V_d(e))=\cl\Int Z(e)$.
\end{proof}

\vspace{.2in}
\section{ $c$-subgroups}

Recall the construction of $c$-subgroups from \cite{lm}.
Given any convex $\ell$-subgroup we let $H_c$ be the convex $\ell$-subgroup generated by the complemented elements in $H$. It follows that $H_c\leq H$ and $H$ is said to be a {\bf $c$-subgroup} if $H=H_c$. Any subgroup containing a weak order unit is a $c$-subgroup and so we consider those $c$-subgroups without a weak order unit as {\bf proper}. Maximal proper $c$-subgroups exist and every proper $c$-subgroup can be extended to one. We denote the set of all maximal $c$-subgroups by $\Max_c(G)$.  An interesting note is that objects in $\Max_c(G)$ need not be prime subgroups. What is also interesting is that the collection of $c$-subgroups forms an algebraic subframe of $\cC(G)$ which we shall by $\cC_c(G)$ (see \cite[Section 3]{lm}).

We wish to view $\Max_c(G)$ as a topological space. Of course, the hull-kernel topology would be the natural one, but we must take care in defining the topology since the maximal proper $c$-subgroups need not be prime subgroups, and so $U_c(g)\cap U_c(h)$ need not equal $U_c(g\wedge h)$ where
$$U_c(g)=\{N\in\Max_c(G): g\notin M\}.$$
The simple twist is to restrict from $G^+$ to $c(G)$.

\begin{proposition}
Let $e,f\in c(G)$, then $U_c(e)\cap U_c(f)=U_c(e \wedge f)$ and $U_c(e)\cup U_c(f)=U_c(e \vee f)$. In particular, for each $M\in\Max_c(G)$ and any complemented pair $e,f\in c(G)$ exactly one of $e$ or $f$ belongs to $M$.
\end{proposition}

\begin{proof}
By convexity, $U_c(e \wedge f)\subseteq U_c(e)\cap U_c(f)$. For the reverse inclusion, let $M\in\Max_c(G)$ and suppose that $e,f\in c(G)$ and $e,f\notin M$. It follows that the $c$-subgroup generated by $e$ and $M$ contains a weak order unit. But this subgroup is $G(e)+M$, where $G(e)$ is the convex $\ell$-subgroup generated by $e$. Similarly, $G(f)+M$ also contains a weak order unit. The intersection of the two subgroups will also contain a weak order unit. However,
$$(G(e)+M)\cap (G(f)+M)= (G(e)\cap G(f)) +M = G(e\wedge f) +M.$$
It follows that $e\wedge f\notin M$.
\vspace{.1in}

The equality involving unions holds in general for any pair of elements of $G^+$.
\end{proof}

\begin{proposition}\label{max_c}
The collection $\{U_c(e):e\in c(G)\}$ forms a base for a topology on $\Max_c(G)$. This topology makes $\Max_c(G)$ into a compact zero-dimensional Hausdorff space.
\end{proposition}

The following result allows us to recognize the relationship between $\Max_d(G)$ and $\Max_c(G)$.
\begin{proposition}
For each $M\in \Max_d(G)$, $M_c \in \Max_c(G)$.
\end{proposition}
\begin{proof}
If $M_c$ is not a maximal proper $c$-subgroup, then it is contained in one, say $T\in \Max_c(G)$, which cannot be contained in $M$ since $M_c$ is the largest $c$-subgroup contained in $M$.  So there must be a complemented element, say $e$, belonging to $T$ which is not in $M$. But then $M$ contains any complement of $e$, which then belongs to $M_c$ and therefore also $T$.  The join of these complemented elements is a weak order unit in $T$, a contradiction.
\end{proof}

We can thus define a map $\rho: \Max_d(G) \ra \Max_c(G)$ by $\rho(M)=M_c$. This map is surjective since any object $N\in\Max_c(G)$ must be contained in a maximal $d$-object since $N$ does not contain any weak order units.

\begin{lemma}\label{cont}
The map $\rho: \Max_d(G)\ra \Max_c(G)$ is continuous.
\end{lemma}

\begin{proof}
Let $U_c(e)$ be a basic open subset of $\Max_c(G)$ with $e\in c(G)$. Let $O=\rho^{-1}(U_c(e))$ and take $M\in O$. Then $M_c\in U_c(e)$ which means that $e\notin M_c$, and also $e\notin M$. So $O\subseteq U_d(e)$. If $e\notin M$, then $e\notin M_c$, showing the reverse inclusion. It follows that $O=U_d(e)$, a basic open subset of $\Max_d(G)$.
\end{proof}

\vspace{.2in}

Here is an alternative way of viewing $\Max_c(G)$; we relate $\sG(G)$ to $\Max_c(G)$.
\vspace{.2in}

Let $H$ be a proper $c$-subgroup and consider
$$\cC[H]=\{\cl\Int Z(e): e\in H^+\cap c(G)\}.$$
Then $\cC[H]$ is a proper filter on the boolean algebra $\sG(G)$. Conversely, if $\cF$ is a filter on $\sG(G)$, then
$$\cC^{-1}[\cF]=\{g\in G: \hbox{ there is some } e\in c(G) \hbox{ such that } \cl\Int Z(e)\in \cF \hbox{ and }g\leq e\}$$
is a proper $c$-subgroup of $G$. These operators are order preserving and so there is a bijection between $\Max_c(G)$ and $\Ult(\sG(G))$; in fact, the bijection is a homeomorphism between the two compact zero-dimensional spaces.

\vspace{.2in}

\section{Main Theorem}

We now come to the main theorem of the article. We would like to mention that our proof avoids the use of the cloz-cover of $YG$ as is done in \cite[Theorem 3.15]{hvw2}. We leave such a construction to another time.

\begin{theorem}\label{max zd}
Let $(G,u)$ be a {\bf W}-object. The following statements are equivalent.
\begin{enumerate}[label={\rm (\arabic*)}, nolistsep]
  \item $\Max_d(G)$ is a zero-dimensional space.
  \item For each pair of distinct maximal $d$-subgroups, say $M_1\neq M_2 \in \Max_d(G)$, there is a complemented element $0<e\in G^+$ such that $e\in M_1\ssm M_2$.
  \item Whenever $k\in G^+$ and $M\in U_d(k)$ there is a $0<e\in c(G)$ such that $e\notin M$ and $0<e\leq k$,
  \item $\rho:\Max_d(G)\ra \Max_c(G)$ is injective.
  \item The collection $\sG(G)$ forms a base for the closed sets on $YG$.
  \item The collection $\{\coz(e): e \in c(G)\}$ is a base for the open sets of $YG$.
\end{enumerate}
\end{theorem}

\begin{proof}
(1) is equivalent to (2). This is obvious as (2) is stating that $\Max_d(G)$ is totally disconnected, which for compact Hausdorff spaces is equivalent to being zero-dimensional.
\medskip

(1) implies (3). Suppose $\Max_d(G)$ is zero-dimensional and let $0<k\in G^+$ and an $M\in\Max_d(G)$ such that $k\notin M$. Without loss of generality, assume that $k$ is not a weak order unit. By hypothesis, there is a complemented element, say $e'\in c(G)$, such that $M\in U_d(e')\subseteq U_d(k)$. Then, $e=e'\wedge k$ is also a complemented element (as $U_d(e)=U_d(e')$) and $e'\notin M$.
\medskip

(3) implies (1). Let $k\in G^+$ and $M\in U_d(k)$. Without loss of generality, assume that $0<k$. Now, $k\notin M$ and so by (3) there is an $0<e\leq k$ with $e\in c(G)$ and $e\notin M$. Now, $M\in U_d(e)\subseteq U_d(k)$, whence $\Max_d(G)$ is zero-dimensional.
\medskip

(1) implies (4). Suppose that $\Max_d(G)$ is zero-dimensional and let $M,N\in \Max_d(G)$ be distinct maximal $d$-subgroups. Then there is some complemented element, say $0<e\in G$ such that $e\in M\ssm N$. It follows that $e\in M_c \ssm N_c$, whence $\rho$ is injective.
\medskip

(4) implies (1). This is the easy case. If $\rho$ is injective, then it is a homeomorphism. But then $\Max_d(G)$ is zero-dimensional since $\Max_c(G)$  is by Proposition \ref{max_c}.
\medskip

(1) implies (5). Suppose that $\Max_d(G)$ is zero-dimensional and let $p\in YG$ and $V\subseteq YG$ a closed subset not containing $p$. By the Yosida Embedding Theorem, there are disjoint $G$-zero-sets, say $Z_1,Z_2$ such that $p\in \Int Z_1$ and $V\subseteq \Int Z_2$. Since $Z_1,Z_2$ are disjoint in $YG$, then $\lambda_u^{-1}(Z_1)$ and $\lambda_u^{-1}(Z_2)$ are disjoint closed subsets of $\Max_d(G)$. By an application of the hypothesis, there is a clopen subset of $\Max_d(G)$, say $B$, such that $\lambda_u^{-1}(Z_1)\subseteq B$ and $\lambda_u^{-1}(Z_2)\cap B=\emptyset$. Set $B'=\Max_d(G)\ssm B$. By Lemma \ref{maxclopen}, $B=U_d(e)$ and $B'=U_d(f)$ where $e$ and $f$ are a complementary pair. Then $\lambda_u(B)\subseteq \cl \coz(e)=\cl\Int Z(f)$ and $\lambda_u(B')=\cl\Int Z(e)$.

Then $p\in \Int Z(f)$ while $V\subseteq \Int Z(e)$. Since $e,f$ are a complementary pair we know that $\Int Z(f)\cap \Int Z(e)=\emptyset$. Consequently, $p\notin \cl\Int Z(e)$ while $V\subseteq \cl\Int Z(e)$.
\medskip

(5) implies (2). Let $M,N\in \Max_d(G)$ be distinct maximal $d$-subgroups. Then there are $0< g,h\in G^+$ such that $g\in M$, $h\in N$ and $\cl \Int Z(g)\cap' \cl\Int Z(h)=\emptyset$. Thus, $\Int Z(g)\cap \Int Z(h)=\emptyset$. Choose $p\in \Int Z(g)$ such that $p\notin Z(h)$. By hypothesis, there is a complemented element $e\in c(G)$ such that $p\notin \cl\Int Z(e)$ and $Z(h)\subseteq \cl\Int Z(e)$. It follows that $e\in N\ssm M$.
\medskip

(5) implies (6). Let $p\in YG$ and $V$ a closed subset of $YG$ not containing $p$. By (5) together with the separation property of $YG$, there is some $e\in c(G)$ such that $V\subseteq \cl\Int Z(e)$. There is some $g\in G^+$ such that $p\notin Z(g)$ and $\cl\Int Z(e)\subseteq Z(g)$. It is straightforward to check that $\cl\Int Z(e+g)=\cl\Int Z(e)$, whence we can assume that $p\notin Z(e)$. It follows that $p\in \coz(e)$ and $\coz(e)\cap V=\emptyset$.
\medskip

(6) implies (5). Suppose the collection $\{\coz(e):e\in c(G)\}$ forms a base for the open subsets of $YG$. Notice that this collection is closed under finite unions and finite intersections. Let $V$ be a closed subset of $YG$ and $p\notin V$. Since $YG$ is compact Hausdorff, and thus $V$ is compact, we can find disjoint members in our base, say $\coz(e)$ and $\coz(f)$ such that $p\in \coz(f)$ and $V\subseteq \coz(e)$. It follows that $p\notin \cl\coz(e)$. Since $e\in c(G)$, $\cl\coz(e)\in \sG(G)$.

\end{proof}

\begin{example}
It is a consequence of condition (4) of Theorem \ref{max zd} and Lemma \ref{cont}, that when $\Max_d(G)$ is zero-dimensional, then so is $\Max_c(G)$. The converse is not true. If $G$ is a {\bf W}-object with the property that the only complemented elements are the weak order units and $0$, then $\Max_c(G)=\{0\}$ which is zero-dimensional. But $\Max_d(G)$ need not be.

For a concrete example, take $X$ to be an almost $P$-space, that is, there are no proper dense cozero-sets, then it follows that the only complemented cozerosets of $X$ are the clopen ones. There exist connected almost $P$-spaces.

\end{example}


\vspace{.2in}











\begin{thebibliography}{99}
\bibitem{banaschewski}
    Banaschewski, B.
    {\it On Wallman's method of compactification.}
    Math. Nachr {\bf 27} (1963) 105--114.


\bibitem{bm4}
    Bhattacharjee, P. and W. Wm. McGovern.
    {\it Maximal $d$-subgroups and ultrafilters}, Rendiconti del Circolo Matematico di Palermo Series, {\bf 67}, (2018) 421--440.


\bibitem{hkm}
    Hager, A. W., C. M. Kimber, and W. Wm. McGovern,
    {\it Least integer closed groups},
    Ordered Algebraic Structures, Dev. Math., {\bf 7}, Kluwer Acad. Publ., Dordrecht, (2002) 245--260.

\bibitem{hkm2} Hager, A. W., C. M. Kimber, and W. Wm. McGovern.
    {\it Weakly least integer closed groups},
    Rend. Circ. Mat. Palermo {\bf 52} (2) (2003) 453--480.

\bibitem{hvw2}
    Henriksen, M., J. Vermeer, and R. G. Woods.
    {\it Wallman covers of compact spaces}.
    Diss. Math. {\bf 280} (1989), 5--34.


\bibitem{lm}
    Lafuente Rodriguez, R. and W. Wm. McGovern.
    {\it When $\Min(G)$ has a clopen $\pi$-base}.
    Mathematics Bohemica, {\bf 146} (2021) 69--89.

\bibitem{wallman}
    Wallman, H.
    {\it Lattices and topological space}.
    Ann. Math. {\bf 39} (1938) 112-126.

\end{thebibliography}
\end{document}